%% file: fantappie.tex
\documentclass[11pt]{article}
\usepackage{latexsym}
\usepackage{amssymb}
\usepackage{euscript}

\let\cal\mathcal
\input{def_latex2}

\newcommand\fant{Fantappi\`e\ transform\ }
\newcommand\fantd{Fantappi\`e\ transform}

\newcommand\al{\alpha}
\renewcommand\O{\Omega}

\newcommand\Opo{\Omega^\circ}
\newcommand\Gp{G^\circ}

\newcommand\Gpl{G_{m+1}}
\newcommand\Gn{G_m}
\newcommand\Gmi{G_{m-1}}
\newcommand\bw{\bar{w}}
\renewcommand\N{{\mathbb N}}
\renewcommand\L{L_0}
\newcommand\sql{\sqrt{L}}
\renewcommand\F{{\cal F}}

\newcommand\ho{{ H^1}}
\newcommand\hop{{ H^1(\Opo)}}
\renewcommand\hi{H^\infty}
\newcommand\vt{\vec{v}}
\newcommand\Op{O^+(\B)}
\newcommand\Sp{S^+(\B)}
\newcommand\Mp{M^+(\B)}
\newcommand\Opd{O^+(\B)^\dagger}
\newcommand\Spd{S^+(\B)^\dagger}
\newcommand\Mpd{M^+(\B)^\dagger}
\renewcommand\LL{{\cal L}}
\newcommand\W{{\mathbb W}}

\begin{document}
\title{Positivity aspects of the \fant}
\author{
John E. M\raise.5ex\hbox{c}Carthy
\thanks{Partially supported by National Science Foundation Grant
DMS 0070639}\\
Washington University\\
St. Louis, Missouri 63130 \\
{\small mccarthy@math.wustl.edu} \and Mihai Putinar
\thanks{
Partially supported by National Science Foundation Grant
DMS 0100367}\\
University of California\\
Santa Barbara, CA 93106 \\
{\small mputinar@math.ucsb.edu} }

\date{March 17, 2004}

\bibliographystyle{plain}

\maketitle
\begin{abstract}
\end{abstract}

\baselineskip = 18pt

\setcounter{section}{-1}
\section{Introduction} Let $\langle z,w\rangle$ be the Hermitian
product between two vectors $z,w \in \C^n, \ n \geq 1$. The \fant
of a complex measure $\mu$, in the terminology of this article, is
the analytic function
$$(\F \mu) (z) \= \int_{\C^n} \frac{d\mu(w)}{1
- \la z,w \ra}, \ \quad \langle z, {\rm supp}(\mu) \rangle \neq 1.$$ For
instance, if ${\rm supp}(\mu)$ is contained in the closure of the
unit ball $\B$, then $(\F \mu) (z)$ is well defined for $z \in
\B$. Note the dimensionless character of the transform, and the
fact that in dimension one ($n=1$) it is essentially 
the Cauchy transform.

The \fant is one of the basic integral operators in the analysis
of several complex variables. It was used for instance in integral
representation formulas for complex analytic functions on convex
domains, in Grothendieck-K\"othe type dualities, in the study of
analytic functionals and in connection with the complex Radon
transform, see \cite{mar63, ai66, henlei84, legr86, ran86, hen90,
hor94}.

On the other hand, in the last decade, the Hilbert function space
supported by $\B$, with reproducing kernel $\frac{1}{1 - \la z,w
\ra}$, was the focus of several investigations contingent to
operator theory, bounded analytic interpolation, factorization,
and realization theories, see \cite{ampi, aem02} and the
references cited there. We call this space the Drury space of the
ball, in honor of the author who first study it \cite{Dru78}. To
put everything into a single sentence, Drury's space turned out to
be universal for multivariate bounded analytic interpolation and
extension results, see \cite{agmc_cnp} for the precise statement.

One of the aims of the present article is to show that there is no
accident that the \fant and the Drury space $\h$ of the ball in
$\C^n$ share the same kernel. Among other observations, we link
the two in a characterization of the images of functions in the
Bergman space of the ball via the \fantd.

The same Hilbert space approach to the \fant gives a conceptually
simple proof of the Martineau-Aizenberg duality theorem: {\it given a
bounded convex domain $\O \subset \C^n$, the \fant establishes a
continuous bijection between the space of analytic germs on
$\overline{\O}$ and the space of analytic functions on the ``dual"
$\O^\circ = \{ z ; \ \langle z, \O \rangle \neq 1 \}$}, see
\cite{mar63, ai66}, and for a variety of other proofs
\cite{legr86, hor94}. The main idea behind our proof is to
consider a relatively compact domain $\O$ in the ball and the
restriction operator $R$ between the Drury space and the Bergman
space $A^2(\O)$ of $\O$. Then $R$ is compact and the
eigenfunctions of its modulus $R^\ast R$ analytically extend to
$\O^\circ$ and diagonalize the \fant on $A^2(\O)$. In such a way
we obtain a familiar picture in the spectral theory of symmetric,
unbounded operators, namely a Gel'fand triple:
$$ \mathcal{F} A^2(\O) \rightarrow \h \rightarrow A^2(\O),$$
where the arrows are restriction maps. Then this duality pattern
carries over to the corresponding Fr\'echet spaces
$\mathcal{O}(\O^\circ) \rightarrow \h \rightarrow
\mathcal{O}(\overline{\O})$.

The starting point for this project was the problem of
characterizing (as much as possible in intrinsic terms) the
Fantappi\`e transforms $M^+(\B)$ of positive measures $\mu$
carried by the closed ball $\overline{\B}$. This is a moment
problem in disguise, and for some good reasons which will appear
below, a simple solution does not seem to exist. In one complex
variable however, a complete answer is given by the Riesz-Herglotz
theorem and any of its many equivalent statements, such as the
spectral theorem for unitary operators.

A natural duality with respect to the Drury space inner product
pairs the convex cone $M^+(\B)$ with the set of all analytic
functions in the ball having non-negative real part, or
equivalently the set of non-negative pluriharmonic functions in
the ball. The latter cone, especially its geometric convexity
features, remains rather mysterious, see \cite {korpuk, aida76}.
The best one can say from our perspective about non-negative
pluriharmonic functions is a positivity criterion found by Pfister
(\cite{pf62}) and Koranyi-Pukansky (see Theorem 5.3 below); or
that they can be regarded, after a polarization in  double the number
of complex variables, as restrictions to an $n$-plane of
non-negative $M$-harmonic functions, for which a much better
understood potential theory exists \cite{kor65, sto94}.

The same can be said, via duality, about the class $M^+(\B)$: a
function $f$  in $M^+(\B)$ can be extended 
to a class of analytic functions 
in double the number of variables
that is
easier to characterize (say in terms of positive definite
sequences). To give some support for the last statement, we remark
that a Bernstein type theorem for Fantappi\`e transforms of
positive measures, in the real analytic sense, was obtained by
Henkin and Shananin \cite{hensh91}. Specifically, the transforms
$$\int_{\B} \frac{d\mu(w)}{1
-\Re \la z,w \ra}$$ of positive measures $\mu$ can be
characterized by Henkin and Shananin's theorem.

Put in equivalent terms, the above difficulty is a reflection of
the
difference between characterizing 
the complex and real Fourier transforms of a
positive measure:
$$ \int_\B e^{-i \langle z, w \rangle} d\mu(w), \ \ \int_\B  e^{-i \Re \langle z, w \rangle}
d\mu(w).$$ Only the latter have the characterization given by
Bochner's theorem \cite{hor94}.

The recent monograph \cite{aps04} by Andersson, Passare and Sigurdsson contains a
thorough treatment of the \fant. We recommend it to the interested
reader, along with the paper \cite{and90}.

Among the duality computations around $M^+(\B)$ we touch the
positive Schur class (so dear in recent times to operator
theorists) and a matrix realization idea, see Sections 6 and 7.

The last part of the article deals with the double layer type and
also Fantappi\`e type potential
$$\int_\B  \Re \frac{1} {1
- \la z,w \ra} d\mu(w), \ \ z \in \B,$$ and an operator valued
extension of it. This gives an estimate of a symmetrized
functional calculus for systems of non-commuting operators, with
Sobolev type bounds on the joint numerical range. Similar
estimates on the numerical range of a single operator were only
recently discovered \cite {de99, putsa}.

We would like to thank Mats Andersson for several helpful comments on
the paper.

{\small \tableofcontents}

\section{Notation and Formulas}
We shall let $\B$ denote the unit ball in $\C^n$. There are three
Hilbert spaces of analytic functions on $\B$ with which we shall primarily
be
concerned: the Drury space $\h$, the Hardy space $H^2$ and the
Bergman space $A^2$. We can define these in terms of their
reproducing kernels (see \eg \cite{ampi} for a description of how
to pass between a Hilbert function space and its reproducing
kernel):
\beq
k_{\h}(z,w) &\=& \frac{1}{1-\la z,w \ra} \\
k_{H^2}(z,w) &\=& \frac{1}{\left[1-\la z,w \ra\right]^n} \\
k_{A^2}(z,w) &\=& \frac{1}{\left[1-\la z,w \ra\right]^{n+1}} .
\eeq

We shall let $\alpha = (\al_1,\dots,\al_n)$ and $\beta$ be
multi-indices in $\N^n$, where as usual $|\al| = \al_1 + \dots + \al_n$
and $\al! = \al_1! \cdots \al_n!$. Because
$$
\frac{1}{\left[1-\la z,w \ra\right]^d}
\= \sum_{\al} \frac{(|\al|+d -1)!}{\al!} z^\al \bw^\al
$$
for every positive integer $d$, we have (with the appropriate
measure normalizations) \se\att
\begin{eqnarray}
\nonumber
\| z^\al \|^2_{\h} &\=& \frac{\al!}{|\al|!} \\
\label{eqa1}
\| z^\al \|^2_{H^2} &\=& \frac{\al!}{(|\al| + n-1)!} \\
\nonumber
\| z^\al \|^2_{A^2} &\=& \frac{\al!}{(|\al| + n)!}
.
\end{eqnarray}
Moreover, as $ 0 < a < b$
implies
$$
\frac{1}{\left[1-\la z,w \ra\right]^b} \ - \
\frac{1}{\left[1-\la z,w \ra\right]^a}
$$
is a positive kernel (as can be seen by expanding
$$
\frac{1}{\left[1-\la z,w \ra\right]^{b-a}}
$$
in a power series and noting that all the coefficients are
positive),
we have
$$
\h\ \subseteq \ H^2 \ \subseteq \ A^2 .
$$
\bs If $f$ is an analytic function, we shall write $f_\al$ for the
coefficient of $z^\al$ in its Taylor expansion at the origin. \bs
The {\it Euler operator} $\L$ is defined by
$$
\L \= \prod_{j=1}^n (j \, + \, \sum_{i=1}^n z_i
\frac{\partial}{\partial z_i} ),
$$
so
$$
\L z^\al \= (|\al| + 1) \dots (|\al| + n) z^\al .
$$
\bs The {\it \fant}  of a complex measure $\mu$
is defined as:
\be
\label{eqa05}
(\F \mu) (z) \= \int_\B \frac{d\mu(w)}{1 - \la z,w \ra} .
\ee
If $f$ is a function in $A^2$, by the \fant of $f$, written $\F f$,
 we mean the
\fant of $f V |_{\B}$, where $V$ is Lebesgue measure in $\C^n$
normalized so that $V(\B) = 1$.
If $f$ is a function in $A^2(\O)$ for some $\O$ in $\B$,
by $\F f$ we mean the \fant of $ f V |_{\O}$ (we shall assume that
the domain $\O$ is clear from the context).
If $f(z) = \sum f_\al z^\al$ is in $A^2$, then
\be
\label{eqfan1}
(\F f )(z)  \= \sum  \frac{f_\a}{(|\a| +1) \dots (|\a|+n)} \, z^\a .
\ee

Thus the \fant operator is positive and compact when restricted to
the Hardy, Bergman or similar function spaces on the ball.

\section{The Euler operator and the \fant}

In this section we identify the inverse of the \fant on the
Bergman space with a partial differential operator, see also
\cite{ai66,hen90} for other formulas of inversion of the \fantd.
The language of unbounded operators is the most appropriate for
our considerations.

The Euler operator $\L$ is a symmetric, positive, densely defined
operator on $A^2$.
Moreover, for $p,q$ polynomials, we have
\se\att\begin{eqnarray}
\nonumber
\la \L p, q \ra_{A^2} &\=& \sum p_\al \bar{q_\al} \frac{\al !}{|\al
|!} \\
&=& \la p, q \ra_\h .
\label{eqeul}
\end{eqnarray}
Let $L$ be the (unique) self-adjoint extension of $\L$ in $A^2$.
\bprop
\label{propdomeul}
The space $\h$ is the domain of $\sql$.
\eprop
\bp
By (\ref{eqeul}), we have
$$
\| \sql z^\al \|_{A^2} \= \| z^\al \|_\h .
$$
So $h \= \sum h_\al z^\al $ is in ${\rm Dom}(\sql)$ iff $\|\sql h
\|_{A^2}$ is finite iff $\| h \|_\h$ is finite. \ep On the other
hand, let $M$ be the self-adjoint extension of $\L$ on $\h$. Then
\beq \la M f, f \ra_\h &\=& \sum (|\a| +1) \dots (|\a|+n) \frac{\a
!}{|\a|!} |f_\a|^2 \\
&=& \sum |(|\a| +1) \dots (|\a|+n) f_\a |^2 \frac{\a!}{(|\a| + n)!}
.
\eeq
So $f$ is in ${\rm Dom}(\sqrt{M})$ iff $f = \F g$ for some $g$ in
$A^2$ (where $g_\al \= (|\a| +1) \dots (|\a|+n) f_\al$).
Summarizing, we have proved:
\bprop
The \fant is a unitary operator from $A^2$ onto ${\rm
Dom}(\sqrt{M})$ in $\h$. We have $M \F = I$ and
\be
\la \F f, g \ra_\h \= \la f, g \ra_{A^2} .
\label{eqfan2}
\ee
\eprop

Thus Drury's space $\h$ in $n$ complex dimensions is a Sobolev
type space of order $n/2$. The \fant is a smoothing operator which
restores, roughly speaking, $n$ radial derivatives of a function
in $A^2$.

\section{The restriction operator}
\label{secc} The restriction operator between two function spaces
is better known for its applications to approximation theory in
one complex variable. We adapt below some basic ideas of
\cite{gps} with the aim of better understanding the \fant on an
arbitrary domain in $\C^n$.

Let $\O$ be a domain whose closure is contained in $\B$, and let
$$
R: \h \ \to \ A^2(\O), \ \ Rf = f|_\O,
$$
be the restriction operator. Because $\O \subset \subset \B$, the
operator $R$ is compact, and so is $R^\ast R : \h \to \h$. Let
$\l_0 = 1 \geq  \l_1 \geq \dots$ be its eigenvalues, and $ 1 \equiv
f_0,f_1, \dots$ be the corresponding eigenvectors, all normalized
to have length $1$.

Note that
\se\att
\begin{eqnarray}
\nonumber
R^\ast R f_k &\= & \l_k f _k \\
\nonumber
\Leftrightarrow\quad \la R^\ast R f_k , g \ra_\h &=& \l_k \la f_k , g
\ra_\h \quad \forall \ g \in \h \\
\label{eqc2}
\Leftrightarrow\quad \la  f_k , g \ra_{A^2(\O)} &=& \l_k \la f_k , g
\ra_\h \quad \forall \ g \in \h.
\end{eqnarray}
In particular, letting $g$ be the reproducing kernel at $z$ for
$\h$, we get \be \label{eqc1} \l_k f_k ( z) \= \int_\O
\frac{f_k(w) dV(w)}{1 - \la z,w \ra} . \ee Equation~\ref{eqc1}
shows in particular that each eigenfunction $f_k$ extends
analytically to the connected component of the origin in $\Opo$,
which is defined by \be \label{eqc17} \Opo \ := \ \{ z \, : \, \la
z , w \ra \neq 1 \ \forall w \in \O \} . \ee Notice too that we
have \be \label{eqc3} \l_k f_k \= \F (  f_k  \, V |_\O). \ee
Define a new Hilbert space $\hop$ of holomorphic functions on
$\Opo$ by asking $\sqrt{\l_k} f_k$ to be an orthonormal basis; so
$$
\hop \= \{ \sum_{k=0}^\i a_k \sqrt{\l_k} f_k \ : \ \sum |a_k|^2 < \i
\}.
$$
By Equation~\ref{eqc2}, the functions $ 1/ \sqrt{\l_k} \, f_k$ are
an orthonormal basis for $A^2(\O)$; by Equation~\ref{eqc3} we see
that
$\F$ maps $1/\sqrt{\l_k} f_k$ to $ \sqrt{\l_k} \, f_k$. Thus we
have:
\bprop
The \fant
is an isometric isomorphism from
$A^2(\O)$ onto $\hop$.
\eprop

Note that the monomials $z^\alpha, \ \alpha \in \N^n$, diagonalize
as before the \fant on a Reinhardt domain $\O$.

In the case of a single complex variable, the eigenfunctions $f_k$
of the modulus of the restriction operator $R^\ast R$ tend to have
some very rigid qualitative properties. They make the core of the
so-called Fisher-Micchelli theory in complex approximation, see e.g.
\cite{gps}.


\section{The dual of $O(\overline{\O})$}
\label{secd}

For this section, fix $G \supset \supset \B$ to be some convex
domain that contains $\overline{\B}$. Let $\O =  \Gp$, defined by
(\ref{eqc17}).
Then $\O \subset \subset  \B \subset  \subset G$.
We want to prove the classical result (see \eg \cite{hor94})
that asserts that the \fant establishes a duality between $O(G)$,
the Fr\'echet space of functions holomorphic on $G$, and
$O(\overline{\O})$, the
space of functions holomorphic on a
neighborhood of $\overline{O}$.

Notice that $\O$ is star-shaped with respect to the origin, but
need not be convex.
\bl
\label{lemd05}
With $G$ and $\O$ as above, $\Opo = G$.
\el
\bp
Fix a non-zero vector $\vt$ in $\C^n$, and consider for what
complex numbers $w$ does $w \vt$ lie in $\O$.
They are precisely the reciprocals of those numbers $z$ such that
$z \vt$ lies in the projection of $G$ onto the (complex) line
through $\vt$.

A point will lie in $\Opo$ if and only if its projection onto every
line lies in the projection of $G$ onto that line. As $G$ is
convex, the Hahn-Banach separation theorem shows that $\Opo = G$.
\ep

Let $\B \subset \subset G_m \subset \subset G_{m+1}$ be a
sequence of smoothly bounded convex domains such that
$$
\bigcup_{m} G_m \= G .\
$$
Let $\O_m = G_m^\circ$; this is a decreasing sequence such that
$$
\bigcap_{m} \O_m \= \overline{\O} .
$$
In Section~\ref{secc} we showed the duality
\beq
A^2(\O_m) \; \times \; \ho(\Opo_m) &\to& \C \\
(f,h) &\=& \la \F f , h \ra_{\ho(\Opo_m)} .
\eeq
For $f$ in $\h$, the pairing is
\beq
(f,h) &\=& \la \F f , h \ra_{\ho(\Opo_m)} \\
&=& \la f, \F^{-1} h \ra_{A^2(\O_m)} \\
&=& \la f, h  \ra_\h .
\eeq
(The last equality can be seen by expanding $f$ and $h$ in terms
of the eigenfunctions $f_k$).
Thus, by general arguments from the theory of
locally convex spaces (see \eg \cite{narb}), there is a
duality pairing
$$
\lim_{\to} A^2(\O_m) \: \times \: \lim_{\leftarrow} H^1 (G_m)
\ \to \ \C .
$$
It is immediate that $$
 \lim_{\to} A^2(\O_m) =
O(\overline{\O}).$$
 To prove that $$
 \lim_{\leftarrow} H^1 (G_m) =
O(G),$$
 we shall use Lemma \ref{lemd1}.
We will let $H^\i(\O)$ denote the space of bounded analytic
functions on $\O$. \bl \label{lemd1} For each $m$, we have the
continuous inclusions
$$
\hi(G_{m+1}) \ \hookrightarrow \ H^1( G_m)
\ \hookrightarrow \ \hi(G_{m-1}) .
$$
\el
\bp
Recall that $\Gmi \subset\subset \Gn \subset\subset \Gpl$.
First we shall prove the inclusion $H^1( \Gn)
\ \hookrightarrow \ \hi(\Gmi)$.

The space $H^1( \Gn)$ is the set
\be\label{eqd1}
\int_{\O_m} \frac{f(w) dV(w)}{1 - \la z,w\ra} ,
\ee
where $f$ ranges over $A^2(\O_m)$.
If $z$ is in $\Gmi$, then the denominator in the integrand is
bounded away from $0$, so
(\ref{eqd1}) gives a function in $\hi(\Gmi)$ whose norm is bounded
by a constant times the norm of $f$ in $A^2(\O_m)$.
\bs
The first inclusion requires more work.
Suppose $h$ is in $\hi(\Gpl)$ and continuous up to the boundary.
Then by the Cauchy integral formula for convex domains
with $C^2$ boundary (see \cite{ran86}[Thm. IV.3.4]) $h$
can be represented as a Cauchy integral of its boundary values:
if $r$ is a defining function for $\partial \Gpl$, then
\be
\label{eqd2}
h(z) \= \frac{1}{(2\pi i)^n} \, \int_{\partial \Gpl}
h(\zeta) \frac{\partial r (\zeta) \wedge (\bar \partial \partial r
(\zeta))^{n-1}}{\la \partial r(\zeta) , \zeta - z \ra^n}.
\ee
The denominator in the integrand is the $n^{\rm th}$ power of the
defining equation for the tangent plane to $\partial \Gpl$ at
$\zeta$. Moreover, $\la \partial r(\zeta) , \zeta \ra$ can never be
zero on $\partial \Gpl$, for otherwise the tangent plane at
$\zeta$ would pass through the origin, contradicting the fact that
$\Gpl$ is convex and $0$ is in its interior.

Therefore, by approximating the integral in (\ref{eqd2}) by a
Riemann sum, we can uniformly approximate $h$
on $\Gn$
by rational functions of the form
\be
\label{eqd25}
\sum_i a_i \frac{1}{(1 - \la z, w^{(i)} \ra )^n}
\ee
where $w^{(i)}$ lie in $\partial \O_{m+1}$ and $\sum |a_i|$ is
uniformly bounded.

Let $k(z,w)$ be the Bergman kernel for $A^2(\O_{m})$.
Then the \fant (on $\O_m$) of $k(z,w)$ is
\be
\label{eqd3}
\int_{\O_m} k(\zeta, w) \frac{1}{1- \la z, \zeta \ra} dV(\zeta)
\= \frac{1}{1 - \la z, w \ra} ,
\ee
and the \fant of any partial differential operator $E$ (in $w$)
applied to $k(z,w)$ is just $E$ applied to the right-hand side of
(\ref{eqd3}).

In particular, let $E$ be the Euler operator adapted to the Hardy
space, {\it viz.}
\be
\label{eqd7}
E_{\bar w} \= \frac{1}{(n-1)!} \
\prod_{j=1}^{n-1} (j \, + \, \sum_{i=1}^n \overline{w_i}
\frac{\partial}{\partial \overline{w_i}} ).
\ee
Then
$$
E_{\bar w} \ \frac{1}{1 - \la z, w \ra} \= \frac{1}{(1 - \la z, w \ra)^n} .
$$
Now, for any $w$ in $\O_{m+1}$, all the functions
$E_{\bar w}\, k(z,w)$
are uniformly bounded in norm in $A^2(\O_m)$.
Therefore (\ref{eqd25}) is the \fant of the
function
$$
\sum a_i E_{\bar w}\, k(z, w^{(i)} )
$$
whose norm is controlled in $A^2(\O_m)$, so
the norm of (\ref{eqd25}) is controlled in $H^1(G_m)$.

Finally, the assumption that $h$ is continuous up to the boundary
of $G_{m+1}$ can be dropped by inserting another smooth
convex domain $G_{m + \frac{1}{2}}$
between $G_m$ and $G_{m+1}$ and
using (\ref{eqd2}) on the boundary of that domain instead.
\ep

We can now prove:
\bt
With notation as above, $\lim\limits_{\leftarrow} H^1 (G_m) =
O(G)$.
\et
\bp
By Lemma~\ref{lemd1}, the spaces
$\lim\limits_{\leftarrow} H^1 (G_m) $ and
$\lim\limits_{\leftarrow} H^\i (G_m)$
are the same.
As a set, $\lim\limits_{\leftarrow} H^\i (G_m)$ consists of those
functions whose restriction to each $G_m$ is bounded and
holomorphic; this is exactly $O(G)$.

The topologies are also the same: by definition, the topology on
$\lim\limits_{\leftarrow} H^\i (G_m)$ is the weakest locally convex
topology such that the restriction maps to every $H^\i (G_m)$ are
all continuous. But this is precisely the topology of uniform
convergence on compact subsets of $G$.
\ep

Thus we have proved:
\bt
\label{thmd2}
With notation as above, the dual of $O(G)$ is
$O(\overline{G^\circ})$. The duality is implemented, for $f$ in
$O(\overline{G^\circ})$ and $h$ in $O(G)$ by choosing $m$
sufficiently large so that $f$ is in $A^2(\O_m)$, and letting
$$
(f,h) \= \la \F f , h \ra_{H^1 (G_m)};
$$
the right-hand side is independent of $m$.
\et

{\bf Remark 1.}
If we approximate $G$ by a decreasing sequence of supersets,
the same method gives a duality between $O(\overline G)$ and
$O(G^\circ)$.

{\bf Remark 2.}
Another formulation of Theorem~\ref{thmd2} is that
$$
O(G) \rightarrow \h \rightarrow O( \overline{\O})
$$
is a {\it Gel'fand triple} (a rigged triple of locally convex spaces)
that diagonalizes the \fantd. See \cite{gv64}. According to the
computations contained in the previous section we also have the
Gel'fand triple
$$
H^1(G) \rightarrow \h \rightarrow A^2( \O) ,
$$
with orthonormal systems $\{ \sqrt{ \l_k } f_k \},\
\{ f_k \},\ \{ 1/\sqrt{\l_k} f_k \}$ and such that
$$
\F ( \frac{1}{\sqrt{\l_k}} f_k) \= \sqrt{\l_k} f_k .
$$

\section{Functions of positive real part}
\label{sece}

We shall use $\Op$ to denote the holomorphic functions of positive
real part on $\B$:
$$
\Op \ :=\ \{ f \, \in \, O(\B) \ : \ \Re(f) \geq 0 \, \}
$$
The following description is due to Kor\'anyi and Pukansky
\cite{korpuk}.

Assume first that $f$ is in $O(\overline{\B})$, and let
$$
S(z,w) \= \frac{1}{(1 - \la z,w \ra)^n}
$$
be the \sz kernel, the kernel for the Hardy space $H^2$.
Then, letting $\sigma$ denote normalized surface area measure on
$\partial \B$, we have
\beq
f(z) S(z,w) &\=& \int_{\partial \B} f(u) S(u,w) S(z,u) d \sigma(u)
\\
\overline{f(w)} S(z,w) &\=& \int_{\partial \B} \overline{f(u)}
S(z,u)S(u,w)d\sigma(u).
\eeq
Hence
$$
[f(z) + \overline{f(w)} ] S(z,w) \= 2 \, \int_{\partial \B}
S(z,u) S(u,w) \Re f(u) d \sigma(u),
$$
or equivalently
$$
\Re f(z) \= \int_{\partial \B}
\frac{S(z,u) S(u,z)}{S(z,z)} \Re f(u) d \sigma(u).
$$
The kernel
$$
P(z,u) \= \frac{|S(z,u)|^2}{S(z,z)}
$$
is called the {\em invariant Poisson kernel} of $\B$, and has been
much studied --- see \cite{kor65, sto94}.

For an arbitrary function $f$ in $\Op$, one considers the dilates
$f_r(z) = f(rz)$ with $r$ increasing to $1$. Then the measures
$\Re f_r \, \sigma$ converge weak-* to a positive measure $\mu$ such that
\be
\label{eqe1}
[f(z) + \overline{f(w)} ] S(z,w) \= 2 \, \int_{\partial \B}
S(z,u) S(u,w) \Re f(u) d \mu(u).
\ee

Notice that the fact that each $\Re f_r (u) $ only has terms in
powers of $u$ and $\overline u$, with no mixed terms, means that
the measure $\mu$ from (\ref{eqe1}) annihilates all monomials
of the form
\be
\left.
\begin{array}{lll}
u^\alpha \overline{u}^\beta  &\quad& \alpha \not\le \beta,\ \beta
\not\le \alpha\\
u^\alpha \overline{u}^{\alpha + \beta} [\alpha_j + \beta_j +1 -
(|\alpha| + |\beta| + n)|u_j|^2 ] &\quad& 1 \leq j \leq n
\\
u^{\alpha+\beta} \overline{u}^{\alpha } [\alpha_j + \beta_j +1 -
(|\alpha| + |\beta| + n)|u_j|^2 ] &\quad& 1 \leq j \leq n .
\end{array}
\right\}
\label{eqe2}
\ee
The first line in (\ref{eqe2}) comes from the fact that there are
no mixed terms in $\Re f_r (u) $, and the second and third
from comparing the
integral $\int_{\partial B} |u^{\alpha+\beta}|^2 d \sigma (u)$
with $
\int_{\partial B} |u^{\alpha+\beta}u_j|^2 d \sigma (u)
$  --- see Formula~\ref{eqa1}.
See \cite{aida76} for a discussion of this point.
We shall call a positive measure on $\partial \B$
that annihilates (\ref{eqe2}) a {\em Kor\'anyi-Pukansky measure}.

We summarize these observations in the following theorem, where
{\em (iii)} is obtained by letting $w=0$ in (\ref{eqe1}).
\bt [Kor\'anyi-Pukansky]
\label{thme1}
Let $f$ be in $O(\B)$. Then the following are equivalent:

(i) The function $f$ is in $\Op$;

(ii) The kernel $[f(z) + \overline{f(w)}] S(z,w) $ is positive
semi-definite;

(iii) There exists a Kor\'anyi-Pukansky measure $\mu$
such that
$$
f(z) \=  \int_{\partial B} [ 2 S(z,u) -1] d\mu(u) \, + \, it
$$
for some $t$ in $\R$.
\et
{\bf Remark 1.} If $\nu$ is an arbitrary positive measure on
$\partial \B$ (\ie not required to annihilate (\ref{eqe2})), then
$$
U(z) \= \int_{\partial B} P(z,u) d\nu(u)
$$
is a non-negative M-harmonic function on $\B$; the converse is
also true --- see \cite{sto94}. Theorem~\ref{thme1} then describes
which such $U$ are also pluriharmonic.
We note that Audibert has a different approach to testing for
pluriharmonicity \cite{au77, rud80}.

\vs
\noindent
{\bf Remark 2.}
Similar results hold, {\em mutatis mutandis,} for the Bergman
kernel.

\section{The Positive Schur Class}
\label{secf}

In this section we shall discuss the {\it positive Schur class of
$\h$}, the class $\Sp$ defined by
$$
\Sp \ := \ \{ f \in O(\B) \, : \,
\frac{f (z) + \overline{f(w)}}{1 - \la z, w \ra} \ {\rm is\
positive\ semidefinite} \}.
$$
These functions are exactly the Cayley transforms of the functions
in what is normally called the Schur class, namely the 
closed unit ball of the multiplier algebra of $\h$.
The fact that $T$ is a contraction if and only if
$(I+T)(I-T)^{-1}$ has positive real part was originally observed by
von Neumann \cite{vonN51}. We shall derive a realization formula
for $\Sp$.

Fix $f$ in $\Sp$. Then there is a Hilbert space $\LL$ and a
holomorphic function $k: \B \to \LL$ such that
\be
\label{eqf1}
\frac{f (z) + \overline{f(w)}}{1 - \la z, w \ra}
\=
\la k(z), k(w) \ra
\ee
(\cite{ampi}[Thm 2.53]).
Particular cases of (\ref{eqf1}) are
\beq
f(z) + \overline{f(0)} &\=& \la k(z), k(0) \ra \\
f(0) + \overline{f(w)} &\=& \la k(0), k(w) \ra \\
f(0) + \overline{f(0)} &\=& \la k(0), k(0) \ra .
\eeq
A little algebraic manipulation gives
$$
\la k(z), k(0) \ra + \la k(0), k(w) \ra - \la k(0), k(0) \ra
\=
(1 - \la z, w \ra) \la k(z), k(w) \ra,
$$
or equivalently
\be
\label{eqf2}
\sum_{i=1}^n \la z_i k(z) , w_i k(w) \ra \=
\la k(z) - k(0), k(w) - k(0) \ra .
\ee
Thus the map
\be
\label{eqf3}
V   \ : \
\left(
\begin{array}{c}
z_1 k(z)\\
\vdots\\
z_n k(z)
\end{array} \right) \ \mapsto \
k(z) - k(0) ,
\ee
defined for $z$ in $\B$, can be extended to an isometry
$$
V \, = \, (V_1, \dots, V_n)  \ : \
\LL^n \to \LL .
$$
The fact that $V$ is an isometry is
expressible as
$$
V_i^\ast V_j \= \delta_{ij} I \quad 1\leq i,j \leq n .
$$
Inverting (\ref{eqf3}), we get
$$
k(z) \= \left(I -  \sum_{i=1}^n z_i V_i \right)^{-1} k(0) .
$$
We have proved
\bt
\label{thmf1}
A function $f$ in $O(\B)$ belongs to $\Sp$ if and only if there
exists an isometry $V \, : \, \LL^n \to \LL$ for some
Hilbert space $\LL$, a vector $\xi$ in $\LL$, and a real number $t$,
 such that
\be
\label{eqf4}
f(z) \=  \la \left[ 2(I - \sum z_i V_i)^{-1} - I \right]
\xi, \xi \ra  \, + \, it.
\ee
\et
{\bf Remark 1.}
The row isometry $V$ in (\ref{eqf4}) is not unique.
It can even be replaced by a row of operators $T = (T_1, \dots,
T_n)$ satisfying
$$
T_i^\ast T_j \quad
\left\{
\begin{array}{lll}
=0, &\quad &i\neq j \\
\leq I &\quad &i = j
\end{array}
\right.
$$
Indeed, following the algebra backwards would give an $f$
satisfying
$$
\frac{f (z) + \overline{f(w)}}{1 - \la z, w \ra}
\ \geq \
\la k(z), k(w) \ra
\ \geq \ 0.
$$

\vs
\noindent
 {\bf Remark 2.} Similar results can be obtained for the
classes
$$
S^+_\alpha (\B) \ := \
\{ f \in O(\B) \, : \,
\frac{f (z) + \overline{f(w)}}{(1 - \la z, w \ra)^\alpha} \ {\rm is\
positive\ semidefinite} \}
$$
for every $0 < \alpha \leq 1$.

\vs
\noindent
 {\bf Remark 3.} Obviously $\Sp \subseteq \Op$. The containment
is proper for $n\geq 2$, as shown by Drury for the inverse Cayley
transforms \cite{Dru78}. Indeed, let
$$
q(z) \= [ \sqrt{n} ]^n \, z_1 \cdots z_n .
$$
Then $|q|$ has supremum $1$ on $\B$, but its multiplier norm
on $\h$ is
$$
\| M_q \| \ \geq \ \frac{\| q \|}{\| 1 \|} \=
\frac{\sqrt{n}^n}{\sqrt{n!}}
$$
by \ref{eqa1}. Therefore
$$
p \= \frac{1+q}{1-q}
$$
is in $\Op$ but not in $\Sp$.

\bs
One can obtain an analogous result to the representation in
Theorem~\ref{thmf1} for $\Op$.
\bt
\label{thmf2}
If a function $f$ belongs to $\Op$ then there
exists an isometry $V \, : \, \LL^{2^{n-1}} \to \LL^{2^{n-1}}$ for some
Hilbert space $\LL$, a vector $\xi$ in $\LL$, and a real number
$t$,
 such that
\be
\label{eqf5}
f(z) \=  \la \left[ 2\left(\ I - \sum_{|\al | \leq\, n, \
|\al| \ {\rm odd}} \sqrt{n \choose |\al|}
\ z^\al V_\al\right)^{-1} - I \right]
\xi, \xi \ra  \, + \, it.
\ee
Here we are writing $V$ as
$$
V \=
\bordermatrix{&\LL& \cdots  &\LL \cr
\LL & & V_\al& \cr
\LL^{2^{n-1} -1} &  & V_\al' &} .
$$
\et
\bp
If $f \in \Op$, then
$$
[f(z) + f(w)^\ast] S(z,w) \= \frac{f(z) + f(w)^\ast}{(1 - \la z,
w\ra)^n} \ \geq \ 0,
$$
and therefore can be factored as $\la k(z),k(w)\ra$
for some $k$ with values in a Hilbert space $\LL$.
We get
$$
\la k(z) - k(0), k(w) - k(0) \ra \=
\left( 1 - (1 - \la z, w \ra)^n \right) \la k(z), k(w) \ra .
$$
Therefore the map 
$$
V \, : \, \bigoplus_{{|\al|\ {\rm odd}} \atop { |\al | \leq n}}
\left( \sqrt{n \choose |\al|} z^\al k(z) \right)
\ \mapsto \
(k(z) - k(0)) 
\bigoplus_{{|\al|\ {\rm even}} \atop{ 0 < |\al | \leq n}}
\left( \sqrt{n \choose |\al|} z^\al k(z) \right)
$$
is an isometry.
Letting $\xi = k(0)$ we get (\ref{eqf5}).
\ep

\section{Herglotz transforms}
\label{secg}

By analogy with Theorems~\ref{thme1} and \ref{thmf1}, we define
$\Mp$ to be the set of functions $f$ that can be represented as
$$
f(z) \= \int_{\partial \B} \left[ \frac{2}{1 - \la z, u \ra} - 1
\right] d\mu(u) \, + \, it
$$
for some positive Borel measure $\mu$ and some real number $t$.
Modulo the harmless imaginary constant, this is just the class of
{\it Herglotz transforms} of positive measures, where the Herglotz
transform of $\mu$ is
$$
(H \mu)(z) \= \int_{\partial B} \frac{1 + \la z, u \ra}{1 - \la z,
u \ra} d\mu(u) .
$$
\bprop
The set $\Mp$ is contained in $\Sp$.
\eprop
\bp
Fix a positive measure $\mu$ on $\partial B$.
Then
$$
\frac{(H\mu)(z) + \overline{(H\mu)(w)}}{1 - \la z, w \ra }
\=
2 \int_{\partial B} \frac{1}{1 - \la z, u \ra} \frac{1 - \la z, u
\ra \la u , w \ra}{1 - \la z, w \ra} \frac{1}{1 - \la u, w \ra} \ d
\mu(u) .
$$
The middle factor is positive semi-definite because for every $u$ in
$\overline{ \B}$, the map
$$
z \mapsto \la z, u \ra
$$
is in the closed unit ball of the multiplier algebra of $\h$.
Therefore the whole expression is
positive, as required.
\ep

Thus we have
$$
\Mp \ \subseteq \ \Sp \ \subseteq \ \Op .
$$
For $n=1$, all three sets are equal, by the Riesz-Herglotz theorem.
For $n > 1$, the second inclusion is strict as was remarked in
Section~\ref{secf}. The main result of this section is that the
first inclusion is also strict.

\bt
\label{thmg1}
For $n \geq 2$, we have $ \Mp \ \subsetneq \ \Sp$.
\et
\bp
By definition, $f$ is in $\Mp$ iff there exists $\mu \geq 0$ such
that
\se\att
\begin{eqnarray}
\label{eqg05}
\frac{1}{2} \left[ f(z) + \overline{f(0)} \right] &\=&
\int_{\partial B} \frac{d\mu(u)}{1 - \la z , u \ra} \\
\label{eqg1}
&=& \sum_{\alpha \in \N^n} z^\alpha
\frac{|\al|!}{\al !} \int_{\partial B} \overline{u}^\al d\mu(u) .
\end{eqnarray}
\att

By Theorem~\ref{thmf1}, $g$ is in $\Sp$ iff there
exists an isometry $V: \LL^n \to \LL$ and $\xi \in \LL$ such that
\se\att
\begin{eqnarray}
\label{eqg15}
\frac{1}{2} \left[ g(z) + \overline{g(0)} \right] &\=&
\la (I - zV)^{-1} \xi, \xi \ra \\
\nonumber
&=& \sum_{j=0}^\infty \la (zV)^j \xi, \xi \ra \\
&=&
\sum_{\alpha \in \N^n} z^\alpha
\la (z^\al)_s (V) \xi, \xi \ra \frac{|\al|!}{\al !} .
\label{eqg2}
\end{eqnarray}
\att
Here
\be
\label{eqg24}
(z^\al)_s (V) \ := \ \frac{\al !}{|\al |!} \sum_i V_{i_1} V_{i_2} \dots
V_{i_{|\al|}}
\ee
is a symmetrized functional calculus, and the $\sum_i$ is the sum
over all permutations of $\alpha_1$ $1$'s, $\alpha_2$ $2$'s, {\it
etc.}
So, for example,
$$
(z_1^2 z_2 )_s (V) \= \frac{1}{3} ( V_1^2 V_2 + V_1 V_2 V_1 + V_2
V_1^2) .
$$

Replacing $\mu$ in (\ref{eqg1}) by $\mu$ composed with complex
conjugation in $\overline{\B}$, we may drop the complex conjugates in the
moments in (\ref{eqg1}), and deduce that $g$ in $\Sp$ lies in $\Mp$
iff there exists a positive measure $\mu$ such that
\be
\la (z^\al)_s (V) \xi, \xi \ra \= \int_{\partial B} z^\al d\mu(z)
\quad \forall \, \alpha \in \N^n .
\label{eqg3}
\ee

In view of M.~Riesz's extension theorem
of positive functionals (see \eg \cite{koo88}), equation~(\ref{eqg3})
is equivalent to
\be
\label{eqg4}
\Re \la p_s (V) \xi, \xi \ra \geq 0
\quad \forall \, p \in \C [z] \cap \Op .
\ee

By Remark 3 after Theorem~\ref{thmf1}, if $S$ is the $n$-tuple of
multiplication by the coordinate functions on $\h$, then
there is a polynomial
$p$ in $\Op$ and a vector $\xi$ for which
$$
\Re \la p_s (S) \xi, \xi \ra
\=
\Re \la p (S) \xi, \xi \ra
\ < \ 0.
$$
By Popescu's dilation theorem for row contractions \cite{po89},
there is a larger Hilbert space $\LL$
containing $\h$ and a row isometry $V$ on $\LL$ such that
$$
\Re \la p_s (V) \xi, \xi \ra \=
\Re \la p_s (S) \xi, \xi \ra
$$
for every $\xi \in \h$.
Therefore (\ref{eqg4}) cannot always hold, and so not all
positive Schur functions are Herglotz transforms.
\ep

\section{Duality results}
\label{sech}

Let us define a sesqui-linear form on $\h$ by
\se\att
\begin{eqnarray}
\nonumber
Q(f,g) &\=& \la f, g \ra_\h \ + \  \overline{f(0)} g(0) \\
&=& \la f + \overline{f(0)}, g \ra_\h
\nonumber\\
&=& \sum f_\alpha \overline{g_\al} \frac{\al!}{|\al|!} + \
\overline{f(0)} g(0).
\label{eqh1}
\end{eqnarray}
When $f$ and $g$ are analytic but
not both in $\h$, we shall use $Q(f,g)$ to
denote the sum (\ref{eqh1}) whenever it converges absolutely;
otherwise we shall consider $Q(f,g)$ undefined.
Given a set ${\cal C}$ in $O(\B)$, we shall let ${\cal C}^\dagger$
denote
$$
{\cal C}^\dagger \, := \, \{ g \in O(B)\, :\, \Re\, Q(f,g)  \geq 0 {\rm \ for \
every\
} f\ {\rm in\ } {\cal C} {\rm \ for\ which\ } Q(f ,g)\ {\rm is\
defined} \}.
$$
\bt
\label{thmh1}
The following dualities hold:

(i) $\Mpd = \Op$;

(ii) $\Opd = \Mp$;

(iii) $\Spd \subseteq \Sp$.
\et
\bp
{\it (i)}
The function $g$ is in $\Mpd$ iff its $\h$-innner product with (\ref{eqg1})
has positive real part for every $\mu$. But this inner product is just $\int_{\partial B} \bar g
d\mu$, so it is neccessary and sufficient that $g$
be in $\Op$.

{\it (ii)}
The function $g$ is in $\Opd$ iff
$$ \Re (f) \geq 0 \ {\rm on\ } \partial \B \quad \Rightarrow \quad
\Re Q(f,g) \geq 0 .
$$
This means that
$$
\Re Q(f,g) = \int_{\partial B} f d\mu
$$
for some positive measure $\mu$.
So $g$ is the Herglotz transform of $(1/2) \mu$.

{\it (iii)} Let $g$ in $\Sp$ be given by (\ref{eqg15}). Then by
(\ref{eqg2}), \be \label{eqh2} Q(f,g) \= 2 \la f_s (V) \xi, \xi
\ra. \ee So $f$ is in $\Spd$ iff (\ref{eqh2}) has positive real
part for every row isometry $V$. By Popescu's theorem again, this
forces $\Re f(S)$ to be positive, where $S$ is the $n$-tuple of
multiplication by the coordinate functions on $\h$. Therefore $f$
must be in $\Sp$. \ep \att \vs 
Theorem~\ref{thmh1} gives another proof that the three classes are
distinct. We do not know if the inclusion
$(iii)$ is proper.
\newline\noindent
{\bf Question {\thetheorem}.}
Is $\Spd = \Sp$?
\vs
Let $E$ be as in (\ref{eqd7}):
\be
\label{eqh4}
E_{z} \= \frac{1}{(n-1)!} \
\prod_{j=1}^{n-1} (j \, + \, \sum_{i=1}^n {z_i}
\frac{\partial}{\partial {z_i}} ).
\ee
Then we can give two further characterizations of the set of Herglotz
transforms of positive measures.

\bt
\label{thmh2}
With notation as above,
\se\att
\begin{eqnarray}
\label{eqh5}
\Mp &\ =\ & \{ f \, : \, \Re \left[ \int (Ef) d \mu \right]
 \geq 0 \  \forall {\rm \ Kor\acute{a}nyi-Pukansky\ }\, \mu \} \\
\nonumber
&=& \{ f \, : \, E_z (f(z)) = v(z,0) \ {\rm for \ some\ non{-}negative}
\att\\
&&\qquad {\rm M{-}harmonic\ } v(z,\bar z) \} .
\label{eqh6}
\end{eqnarray}
\et
\bp
(\ref{eqh5}):
When calculating $Q$, there is no loss of generality in assuming
$f(0)$ and $g(0)$ are real, and we shall assume this below.
Let $g$ be in $\Op$; by
Theorem~\ref{thme1}, there is some
Kor\'anyi-Pukansky measure $\mu$ such that
\be
\label{eqh7}
g(z) + g(0) \= \int_{\partial B} \frac{d\mu(u)}{(1 - \la z, u \ra)^n}.
\ee
So $f$ is in $\Mp = \Opd$ iff
the $\h$ inner product of $f$ with (\ref{eqh7}) always has positive
real part (here the inner product is evaluated formally on power
series).
But, assuming $f$ is regular enough,
\beq
\la f(z) ,
\int_{\partial B} \frac{d\mu(u)}{(1 - \la z, u \ra)^n} \ra_\h &\=&
\la f(z) , E_z
\int_{\partial B} \frac{d\mu(u)}{1 - \la z, u \ra } \ra_\h \\
&\=&
\int_{\partial B} \la E_z f(z), \frac{1}{1 - \la z, u \ra } \ra_\h
d\mu(u) \\
&\=&
\int_{\partial B} (Ef)(u) d\mu(u).
\eeq
Therefore (\ref{eqh5}) holds.

\vs
(\ref{eqh6}):
Let
\be
\label{eqh73}
v(z, \bar z) \ := \ \int \frac{S(z,u) S(u,z)}{S(z,z)} d\mu .
\ee
Then if $f$ is given by (\ref{eqg05}), we find
$$
E_z [ f(z) + \overline{f(0)} ] \= 2
v(z,0) .
$$
Conversely, given $v$, let $\mu$ be defined by (\ref{eqh73}).
Then $f$ given by (\ref{eqg05}) satisifes
(\ref{eqh6}).

\ep
\vs
\noindent
{\bf Example.} By direct calculation,
\be
\label{eqh8}
E_z z^\alpha \=
{n + |\al| -1 \choose n-1} z^\al .
\ee
So if $f$ is the sum of a function $f_d$ homogeneousof degree $d$
and a constant $f_0$, then
$$
E_z f(z) \= f_0 + {n + |d| -1 \choose n-1} f_d (z) \= f(rz),
$$
where
$$
r = {n + |d| -1 \choose n-1}^{1/d} .
$$
A sufficient condition for
(\ref{eqh5}) to hold is that
$$
\Re f |_{B(0,r)} \geq 0 .
$$
\bs
As $\Mp \subsetneq \Sp \subsetneq \Op$, any function $f$ in $\Mp$ must
have realizations as in Theorems~\ref{thmf1} and \ref{thme1}.
How are they related?

Assume $f$ is in $\Mp$, and $\Im f(0) = 0$ for convenience.
Then $f$ is the Herglotz transform of some measure $\mu$ on $\partial
\B$. Therefore, if $N = (N_1, \dots, N_n)$ is the normal $n$-tuple of multiplication by
the coordinate functions in $L^2(\mu)$, we have
\be
\label{eqh9}
f(z)  \=  \la [ 2 (I - z \cdot N)^{-1} - I ] 1, 1 \ra.
\ee
As $N N^\ast \= \sum N_i N_i^\ast = I$, $N$ is a co-isometry, and
therefore a row contraction. Let $V$ be Popescu's 
row isometric dilation of $N$ \cite{po89}.
Then replacing $N$ by $V$ in (\ref{eqh9}) we get the same function;
this is the realization of $f$ in Theorem~\ref{thmf1}.

The connection with $\Op$ is less clear.
If $g$ is a matrix or scalar valued analytic function on $\B$ with 
$\Re g \geq 0$, then
the weak-* limit of the measures $ \Re g(rz) d\sigma(z)$ is a positive
operator valued Kor\'anyi-Pukansky measure $E$ on $\partial \B$
satisfying
$g(z) \= \int 2 S(z,u) -1 dE(u) $.
By Naimark's dilation theorem \cite{fil70}, $E$ has a dilation to a
spectral measure, and so if $N$ is the normal $n$-tuple corresponding
to this spectral measure, and $P$ is the projection onto the range of
$E$, we have
\be
\label{eqh10}
g(z) \= P 2(I-z \cdot N)^{-n} -1 P .
\ee
Notice that if $T$ is any row contraction, then
$$
2 \Re ( I - z \cdot T)^{-1} - I \ \geq \ 0.
$$
So if we let $g(z) = 2(I - z \cdot T)^{-1} - 1$, then 
(\ref{eqh10}) applied to $g(z) + g(0)^\ast$ gives
$$
(I - z \cdot T)^{-1} = P (I - z \cdot N)^{-1} P.
$$

\section{Functional calculus on the numerical range: The Ball}
\label{seci}

Let $R$ be an operator on a Hilbert space $\LL$. Its {\it
numerical range}, denoted $W(R)$, is the set
$$
W(R) \ := \ \{ \la R \xi, \xi \ra \ : \ \| \xi \| = 1 \} .
$$

By a classical Theorem of Hausdorff and Toeplitz, the numerical
range of an operator is a convex set that contains in its closure
the spectrum. Recently, B. and F. Delyon proved that for any $R$,
its numerical range is an $M$-spectral set for some $M$
\cite{de99}, \ie \be \label{eqi1} \| p(R) \| \ \leq \ M \| p
\|_{W(R)} \quad \forall \ {\rm polynomials\ } p . \ee For an
alternative proof, with an analysis of the best $M$, see
\cite{putsa}. In this section, we shall extend (\ref{eqi1}) to
$n$-tuples with numerical range in the ball.

Let $T = (T_1, \dots, T_n)$ be an $n$-tuple of operators on a
Hilbert space $\LL$. {\em We do not assume that the operators
commute with each other.} We shall say the numerical range of $T$
is contained in $\overline{\B}$, written $\W(T) \subseteq
\overline{\B}$, if for every $u$ in $\overline{\B}$,
\be
\label{eqi2}
W(u  T) \= W( \bar u_1 T_1 + \dots + \bar u_n T_n) \ \subseteq \
\overline{\D} .
\ee
Our standing
assumption throughout this section will be that $\W(T) \subseteq
\overline{\B}$.

\bl
\label{lemi1}
$$
\Re (I - \bar u T)^{-1}  \ \geq \ 0 \quad \forall \, u \, \in \, \B
.
$$
\el

\bp
Note that
$$
(I - \bar u T)^{-1} + (I -  u T^\ast)^{-1}
\=
(I - \bar u T)^{-1} \left[ 2 - \bar u T - u T^\ast \right] (I -  u
T^\ast)^{-1} .
$$
The quantity in brackets is positive iff $\Re \bar u T \leq I$;
this holds for all $u$ iff $\W(T) \subseteq
\overline{\B}$.
\ep

Consider the measures $\Re(I - r \bar u T)^{-1} d \sigma(u)$.
These
are all positive by Lemma~\ref{lemi1}, and have total mass $1$.
Therefore the positive operator valued measure
\be
\label{eqi3}
d\mu_T (u) \ := \ {\rm weak}^*-\lim_{r \nearrow 1} \Re
(I - r \bar u T)^{-1} d \sigma(u)
\ee
exists and is well-defined.
Define
\beq
\Xi \, : \, \C[z] &\ \to \ &  B(\LL) \\
p &\mapsto & \int_{\partial B} p \, d\mu_T .
\eeq
If $\Re p \geq 0 $ on $\partial B$,
then $\Re [ \Xi (p)] \geq 0 $.

To understand $\Xi$ better, let us consider the scalar case.
Define
\beq
\La p (z) &\ := \ &
\int_{\partial B} \frac{1}{2} \left[ \frac{1}{1 - \la u, z \ra}
+ \frac{1}{1 - \la z, u \ra} \right] p(u) d\sigma (u)
\\
&=& \frac{1}{2} \F (p d\sigma) \, + \, \frac{1}{2} p(0) .
\eeq
Then
$$
\La\, :\, z^\al \mapsto
\left\{
\begin{array}{ll}
\frac{1}{2} \frac{1}{(|\al | + 1) \cdots  (|\al | + n - 1)}
z^\al \quad & \al \neq 0 \\
&\\
1 & \al = 0 .
\end{array}
\right.
$$
Note that
\be\label{eqi4}
\Xi (p) \= \left[ \La p \right]_s (T) ,
\ee
where the subscript means the symmetrized functional calculus
from (\ref{eqg24}).

By direct calculation on monomials,
the operator
\be\label{eqi5}
\Gamma \, : \ p \ \ \mapsto \  2 (n-1)! E_z [ p - p(0)] \, + \, p(0)
\ee
is the inverse to $\La$.

By Naimark's dilation theorem \cite{fil70}, the positive operator valued measure
$\mu_T$ has a dilation to a spectral measure on $\partial B$,
whose values are projections in  a Hilbert space
$\K \supseteq \LL$. If $P$ is projection from $\K$ onto $\LL$, then
$$
\int p d\mu_T \= P p(N) P
$$
where $N$ is the normal $n$-tuple
of multiplication by the coordinate functions.

By (\ref{eqi4}), for any polynomial $p$,
$$
\Xi (\Gamma \,p ) \= \left[ \La \, \Gamma \, p \right]_s (T) \= p_s(T)
.$$
Therefore
\be
\label{eqi6}
p_s(T) \= \int \Gamma\, p\  d\mu_T \= P \Gamma p (N) P.
\ee
The above argument goes through unchanged if $p$ is matrix-valued.
If $p(z) \= \sum A_\alpha z^\alpha$, with the coefficients $A_\alpha$
$d \times d$ matrices, then 
$p_s(T)$ is the $d \times d$ operator valued matrix with $(i,j)$ entry
$[p_{ij}]_s (T)$, and $\Gamma p$ is likewise obtained by applying
$\Gamma$ entrywise.

We have proved
\bt
\label{thmi1}
Let $T$ have $\W(T) \subseteq \overline{\B}$, and let $\Gamma$ be defined by
(\ref{eqi5}).
Then, for any polynomial $p$, scalar or matrix valued, we have
$$
\| p_s (T) \| \ \leq \ \| \Gamma p \|_{\overline{\B}} .
$$
\et

\noindent
{\bf Example 1.} Let $p$ be homogeneous of degree $d$. Then
$$
\Gamma p \= 2 (d+1) \dots (d+n-1) p ,
$$
so
$$
\| p_s(T) \| \ \leq \ 2 (d+1) \dots (d+n-1) \| p \|_{\overline{\B}}.
$$

\bs
\noindent
{\bf Example 2.}
In one complex variable $(n=1)$ Theorem~\ref{thmi1} gives
$$ \| p(T)\| \ \leq\  \|2 p - p(0) \|_{\overline{\D}} \ \leq \ 3 \| p
\|_{\overline{\D}}$$
whenever $W(T) \subseteq \overline{\D}$.
This was first obtained in \cite{de99} and \cite{putsa}.
Note that the last inequality says that $T$ is completely 
polynomially bounded, since there is one constant that works for all matrix
valued polynomials. Therefore by Paulsen's theorem \cite{pau84} any operator with numerical range in the
closed unit disk is similar to a contraction.
\bs
In light of (\ref{eqi6}), it is natural to ask what conditions on $p$
make $\Re( \Gamma p) $ non-negative. As $\Lambda$ is the inverse of
$\Gamma$, this is equivalent to asking when $p$ is in $\Lambda ( \Op)$.

Now, $q$ is in $\Op$ iff $H[\Re q \, d\sigma]$ is in $\Mp$.
By comparing the formulas for $H$ and $\Lambda$, we get
$$
2 \Lambda[q] \=
H[\Re q \, d\sigma]
\, + \,
i\Im q(0)
\, + \,
q(0).
$$
Therefore $p$ is $\Lambda (q)$ for some $q$ in $\Op$ if and only if
\be
\label{eqi7}
p \, - \, \frac{1}{2} \Re p(0) \, - \, 
i \Im p(0) \ \in \ \Mp .
\ee

\begin{cor}
Suppose 
$\W(T) \subseteq
\overline{\B}$
and $p$\ satisfies (\ref{eqi7}).
Then $\Re p_s(T) \geq 0$.
\end{cor}

\bs
Alternatively, one can work directly with $\Gamma p$, which is easily
calculated when $p$ is decomposed into homogeneous pieces.

\begin{cor}
Let $p = p_0 + p_1 + \dots + p_d$ be the decomposition of $p$ into
homoegeneous polynomials. If 
$$
\frac 12 \, \Gamma p \= \frac 12 p_0 + \frac{n!}{1!} p_1 + \dots
+ \frac{(n+d-1)!}{d!} p_d $$
has positive real part on $\B$, then
$\Re p_s(T) \geq 0 $ for every $T$ with $\W(T) \subseteq
\overline{\B}$.
\end{cor}

\noindent
{\bf Example 3.}
Suppose $T_2 = T_1^\ast$, and the numerical radius of $T$ is at most
$1 / \sqrt{2}$ (the
numerical radius 
is the supremum of the moduli of the numbers in the numerical range).
Then $\W(T_1, T_2) \subseteq \overline{\B_2}$.
Let $n=2$, and, for $m \geq 2$, let 
$$
p(z_1,z_2) \= \frac{2m+1}{2^{m-1}} \, - \, z_1^m z_2^m
. $$
Then $$
\Gamma p\, (z_1,z_2) \= 
\frac{2m+1}{2^{m-1}} \, - \, 2(2m+1) \,z_1^m z_2^m,
$$
which is nonnegative on $\B_2$. 
Therefore $(z_1^m z_2^m)_s (T)$, which is the average of all
${2m} \choose{m}$ ways of writing $m$ copies of $T_1$ and $m$ copies of
$T_1^\ast$, is less than or equal to $(2m+1)/2^{m-1}$ times the
identity:
\be
\label{eqi8}
(z_1^m z_2^m)_s (T_1, T_1^\ast) \ \leq \
\frac{2m+1}{2^{m-1}} I
\ee 
When $m=1$, inequality~(\ref{eqi8}) is worse than the trivial 
one obtained from
the observation that $\|T_1\|$ is at most twice the numerical radius.
For $m \geq 2$, the inequality is better than this trivial one.

If $\dis T_1 \= \left( \begin{array}{cc} 0 & \sqrt{2} \\
0 & 0 \end{array} \right)$, then only two terms in the whole sum
are non-zero. One obtains the inequality
$$
{{2m} \choose m} \ \geq \ \frac{2^{2m-1}}{2m+1} .
$$

\bs
\noindent
{\bf Example 4.} Let $T = (T_1,\dots, T_n)$ be a {\em commuting}
$n$-tuple, with $\W(T) \subseteq \overline{\B_n}$.
Assume also that the operators are jointly nilpotent, in the sense
that
there is some integer $N$ such that $T^\alpha =0$ whenever $|\al| > N $.
Then there is some constant $C_N$ such that if
$p(z) \= \sum A_\alpha z^\al$ is a matrix valued polynomial, then
$$
\sum_{|\al| \leq N} \|A_\al\| \ \leq \ C_N \| p \|_{\B_n} .
$$
Therefore
$$
\| \Gamma p \| \ \leq \ 
2(N+1) \cdots (N+n-1) \, C_N \| p \|_{\B_n}.
$$
Therefore the $n$-tuple $T$ is completely polynomially bounded, and so
it is similar to an $n$-tuple that has a normal dilation to $\partial
\B_n$ by \cite[Thm xx]{pau86}.

\section{General convex sets}
\label{secj}

Let $\O$ be a closed convex set in $\C^n$ with $0$ in the
interior. Assume that $\O$ is {\em circular}, \ie whenever $z$ is
in $\O$, then so is $e^{i\theta} z$ for all $\theta \, \in \,
[0,2\pi]$.
Let $\Opo$ be defined by (\ref{eqc17}). Let $T = (T_1, \dots, T_n)$
be an $n$-tuple of not necessarily commuting operators. Then we say
$\W(T) \subseteq \O$
if
$$
W( \bar u_1 T_1 + \dots \bar u_n T_n ) \ \subseteq \ \overline{\D}
\qquad \forall \, u \, \in \, \Opo.
$$

Fix some $n$-tuple $T$ with $\W(T) \subseteq \O$.
Let $\omega$ be harmonic measure for $\Opo$ at $0$.
As in Section~\ref{seci}, we can define the positive operator
valued measure $\mu_T$
on $\partial \Opo$ by
$$
d \mu_T (u) \ :=
\ {\rm weak}^*-\lim_{r \nearrow 1} \Re
(I - r \bar u T)^{-1} d \omega(u)  .
$$
Define $\Xi$  and $\Lambda$ by
\beq
\Xi (p) &\=& \int p \, d \mu_T
\\
\Lambda p &=& \frac{1}{2} \F( p\, \omega)  + \frac{1}{2} p(0) .
\eeq

If $\Opo$ is Reinhardt (\ie invariant under rotation of each
coordinate separately) then the monomials are orthogonal, and
$$
\F (u^\al \omega) (z) \= \left[ \sqrt{\frac{|\al|!}{\al!}} \int_{\partial
\Opo} |u^\al|^2 d\omega(u) \right] z^\al .
$$
So $\Gamma = \Lambda^{-1}$ exists and is diagonalized by the
orthogonal basis of the monomials.

Even if $\Opo$ is not Reinhardt, it is invariant under the action
of the circle group by hypothesis. Therefore in $L^2(\omega)$,
homogeneous polynomials of different degrees are orthogonal.
Let ${\cal P}_d$ denote the homogeneous polynomials of degree $d$.
Then $\Lambda\, : \, {\cal P}_d \to {\cal P}_d$.
Moreover, $\Lambda$ has no kernel, because
$ \Lambda p \= 0 $
implies $p$ is orthogonal
to every power of $u$, and so $\int |p|^2 d \omega = 0$.
Therefore, as ${\cal P}_d$ is finite dimensional, $\Gamma =
\Lambda^{-1}$ exists and maps ${\cal P}_d$ onto ${\cal P}_d$.
We can therefore repeat the argument of Theorem~\ref{thmi1}, and
get:
\bt
\label{thmj1}
With notation as above, let $T$ have $\W(T) \subseteq \overline{\O}$.
Then, for any polynomial $p$, we have
$$
\| p_s (T) \| \ \leq \ \| \Gamma p \|_{\overline{\Opo}} .
$$
\et

\bibliography{references}

\end{document}

%% file: def_latex2.tex
\usepackage{latexsym}
\usepackage{amssymb}
\usepackage{euscript}
\usepackage[dvips]{graphics}
\usepackage{epsf}

\let\cal=\mathcal      

\def\mcc{M\raise.5ex\hbox{c}C}
\def\mccarthy{M\raise.5ex\hbox{c}Carthy}
\def\Hu{\H}
\def\sz{Szeg\Hu{o} }

\def\eg{{\it e.g. }}
\def\ie{{\it i.e. }}

\def\h{{\cal H}}

\def\K{{\cal K}}

\def\N{{\cal N}}
 
\def\LL{{\cal L}}

\def\hi{H^\infty}

\def\La{\Lambda}

\def\a{\alpha}
\def\l{\lambda}

\let\i=\infty

\def\la{\langle}
\def\ra{\rangle}
\def\={\ = \ }

\def\F{{\cal F}}

\def\C{\mathbb C}
\def\R{\mathbb R}
\def\T{\mathbb T}
\def\D{\mathbb D}
\def\B{\mathbb B}

\def\dis{\displaystyle}

\def\be{\setcounter{equation}{\value{theorem}} \begin{equation}}
\def\ee{\end{equation} \addtocounter{theorem}{1}}
\def\beq{\begin{eqnarray*}}
\def\eeq{\end{eqnarray*}}
\def\se{\setcounter{equation}{\value{theorem}}} 
\def\att{\addtocounter{theorem}{1}}
\def\vs{\vskip 5pt}
\def\bs{\vskip 12pt}

\def\bp{{\sc Proof: }}
\def\ep{{}{\hfill $\Box$} \vskip 5pt \par}

\def\bl{\begin{lemma}}
\def\el{\end{lemma}}
\def\bt{\begin{theorem}}
\def\et{\end{theorem}}
\def\bprop{\begin{prop}}
\def\eprop{\end{prop}}
\def\bd{\begin{definition}}
\def\ed{\end{definition}}
\def\br{\begin{remark}}
\def\er{\end{remark}}
\def\bexer{\begin{exercise}}
\def\eexer{\end{exercise}}

\newtheorem{theorem}{Theorem}[section]
\newtheorem{prop}[theorem]{Proposition}
\newtheorem{lemma}[theorem]{Lemma}
\newtheorem{cor}[theorem]{Corollary}

\newtheorem{definition}[theorem]{Definition}